\documentclass[11pt]{amsart}

\makeatletter
\def\@normalsize{\@setsize\normalsize{10pt}\xpt\@xpt
\abovedlayskip 10pt plus2pt minus5pt\belowdisplayskip
\abovedisplayskip \abovedisplayshortskip \z@
plus3pt\belowdisplayshortskip 6pt plus3pt
minus3pt\let\@listi\@listI}
\def\subsize{\@setsize\subsize{12pt}\xipt\@xipt}
\def\section{\@startsection {section}{1}{\z@}{1.0ex plus 1ex minus
2ex}{.2ex plus .2ex}{\large\bf}}
\def\subsection{\@startsection {subsection}{2}{\z@}{.2ex plus 1ex}
{.2ex plus .2ex}{\subsize\bf}} \makeatother

\newtheorem{theorem}{Theorem}

\newtheorem{corollary}{Corollary}

\newtheorem{example}{Example}

\newtheorem{proposition}{Proposition}

\theoremstyle{remark}

\begin{document}

\baselineskip=12pt

\address{
UniBW, Munich, Germany }

\email{Alexey.Ostrovskiy@UniBw.de}
\subjclass[2000]{Primary  28A05, 54C10, 54E40, 26A21  Secondary 54C08  }

\keywords{}

\title{\bf {Closed-Constructible functions are  Piece-Wise Closed  
}}

\author {Alexey Ostrovsky}
\maketitle

\begin{abstract}

A subset    $B   \subset Y$   is  constructible if it   is
an element of the smallest family that contains all open sets and
is stable under finite intersections and complements. A function
$f : X \to  Y$ is said to be piece-wise closed  if  $X$  can be
written as a countable union of closed sets  $Z_n$ such that $f$
is closed on every $Z_n.$  We prove that if  a  continuous
function $f$   takes
each closed set into a constructible subset of $Y$, then $f$  is  piece-wise
closed. 

 


\end{abstract}

   \vspace{0.15in }
  All spaces in this paper are supposed to be separable and metrizable, and all the functions are supposed to be continuous and
onto.

          \vspace{0.15in }

 A subset    $B   \subset Y$  is  \emph{constructible} if it   is  an element of the smallest family that contains all open sets and is stable under
 finite intersections and complements.

In topology, a constructible set  is a finite union of locally
closed sets (a set is  \emph{locally closed}  or is an
\emph{$LC_1$}-set if it is the intersection of an open set and a
closed set).

In algebraic geometry, a constructible set is any zero set of a system of polynomial equations
and inequations.
A function $f$  is said to be \emph{closed--constructible}  (resp.,  \emph{closed--$F_{\sigma}$}) if  $f$ takes   closed  sets into constructible
 (resp., $F_{\sigma}$)  ones.
 It is clear that every  closed--constructible function  is  closed--$F_{\sigma}.$

  A function $f : X \to  Y$ is said to be  \emph{ piece-wise closed}
  if $X$ can be written as a countable union of closed sets $Z_n$ such that $f$ is
closed on every $Z_n.$

  \vspace{0.05in}

  Hansell,      Rogers and Jayne      gave  a corrected form of  their previous  result 
        \cite[Theorem 1]{R},  \cite{R2} using additional   hypotheses   $(a)-(d)$    \cite[Theorem 3]{H}.

        \vspace{0.05in}

        Under hypothesis  $(a)$ (=Fleissner's axiom, which is consistent with the usual axiom ZFC), they established the correctness of  their first conclusion \cite{R}[ Lemma 2]:
      each continuous, closed-$F_\sigma$ function between absolute Souslin sets $X$ and $Y$  is piece-wise closed.

      \vspace{0.15in }

       Unfortunately, the theory above is not sufficient for  important  applications such as the case of  closed-$LC_1$ functions between non-Souslin  subsets of the real line $\textbf{R}.$
        \vspace{0.05in }

        Motivated by this
observation, we study the extensions of the theory of
closed--Borel functions, whose major case is closed-$LC_2$
functions (a  simple case for such non-continuous functions  was
recently considered in \cite{2}).

 We    obtain the following main  theorem: 


  \begin{theorem}
  Every   closed--constructible  function    is  piece-wise closed.

  \end{theorem}

   Theorem 1 immediataely implies the following:

  \begin{corollary}

 Each   constructible-measurable,  closed  (or open),  one-to-one  function  is  piece-wise  continuous.

    \end{corollary}

Note that this result is  different from that of  Banakh and
Bokalo \cite{TB}, which states that a constructible-measurable
function of hereditary Baire space $X$ is piece-wise  continuous.

  \vspace{0.15in}

   Finally, we give some simple observations  to   the   hypothesis  $(d)$ mentioned   in the beginning:
   each  preimage $f^{-1}(y)$ of points from $Y$ is compact.

         \vspace{0.10in}

         This hypothesis
            looks   fairly strong: as we demonstrated in simple  Proposition 3, a function $f$ under hypothesis $(d)$ becomes
            a closed--$F_{\sigma}$ function.  However, in the same situation     (Example
            3), such a conclusion  becomes false if the requirement $(d)$ is weakened to the following:  each  preimage $f^{-1}(y)$ of points from  $Y$ is completely metrizable.

 \begin{proposition}    Let
$f:
X  \to  Y$ be a continuous  function with compact fibers and
 $f$ take    elements of a clopen base $ \mathcal{B}$ into  closed sets.  Then $f$  is a  closed and, hence,  closed--constructible function.
 \end{proposition}

Note that in following Example 1  all   preimages $f^{-1}(y)$ of points from $Y$    are completely metrizable, but   not compete under the given metric.  

 \begin{example}

 There exists a continuous  and  open function  $f: I_X  \to I_Y  $ from Polish spaces  $ I_X ,I_Y \subset  \textbf{C}$ that takes  elements  of some clopen base  of $I_X$  into clopen  sets, but $f$ is not   closed--$F_{ \sigma}.$ 

\end{example}

\end{document}